\documentclass[12pt]{article} 
\usepackage{latexsym, 
amsmath, 
amscd, 
amsthm, 
graphicx, epsfig, amssymb}

\newtheorem{thm}{Theorem}[section]

\newtheorem{cor}[thm]{Corollary}

\newtheorem{rema}[thm]{Remark}

\newtheorem{conj}[thm]{Conjecture}

\newtheorem{expl}[thm]{Example}

\newcommand{\C}{\mathbb{C}}

\newcommand{\R}{\mathbb{R}}

\newcommand{\N}{\mathbb{N}}

\newcommand{\one}{\mathbf{1}}

\title{ {\bf Representations of vertex operator algebras and braided
finite tensor categories} }
\date{}
\author{Yi-Zhi Huang}
\begin{document}

\bibliographystyle{alpha}
\maketitle

\begin{abstract}
We discuss what has been achieved in the past twenty years 
on the construction and study of a
braided finite tensor category structure on a suitable module category
for a suitable vertex operator algebra.
We identify the main difficult parts
in the construction, discuss the methods developed to 
overcome these difficulties and 
present some further problems that still need to be solved.
We also choose to discuss three among the numerous 
applications of the construction.  
\end{abstract}

\renewcommand{\theequation}{\thesection.\arabic{equation}}
\renewcommand{\thethm}{\thesection.\arabic{thm}}
\setcounter{equation}{0}
\setcounter{thm}{0}

\section{Introduction}

Finite tensor categories, roughly speaking,
are rigid tensor categories satisfying 
all reasonable finiteness conditions. The category of finite-dimensional 
representations in positive characteristic of a finite group 
is an example of a finite tensor category. Such a finite 
tensor category is symmetric. Another 
class of examples of finite 
tensor categories is constructed from 
representations of  quantum groups
at roots of unity.  In this case, the finite tensor categories are
not symmetric but are instead braided. Moreover, 
some of them are modular tensor categories \cite{T1} \cite{T2}
which satisfy additional conditions,
including, in particular, semisimplicity (of modules) and a nondegeneracy 
property. 
Finite tensor categories have been studied systematically by 
Etingof and Ostrik \cite{EO}. In general, finite tensor categories
are not necessarily semisimple. 

In the semisimple case, modular tensor categories play 
an important role in 
the study of quantum groups, knot and three-manifold invariants,
three-dimensional topological quantum field theories and 
rational conformal field theories. They arose first  
in the study of rational conformal field theories. 
In \cite{MS1} \cite{MS2}, Moore and Seiberg derived 
a set of polynomial equations from an axiom system for a rational 
conformal field theory. Moreover, after Witten commented that one
of their equations was analogous to Mac Lane's coherence property,
they demonstrated  in these papers a convincing analogy between
the theory of such polynomial equations and 
the theory of tensor categories.
Later the mathematical notion of modular tensor category 
based on the theory of tensor categories
was formulated precisely by Turaev in \cite{T1} and \cite{T2}. 
Examples of modular tensor categories were constructed from 
representations of quantum groups but the problem of constructing 
modular tensor categories from candidates for conformal field theories,
especially the proofs of the rigidity and the 
nondegeneracy property, 
was open for a long time.  It was solved by the author in 2005
in \cite{H11} (see also the announcement \cite{H8} and the expositions
\cite{H9} and \cite{Le})
using the representation theory of vertex operator algebras,
including, in particular, the semisimple tensor product theory 
developed by Lepowsky and the author \cite{tensor0}--\cite{tensor3}
\cite{tensor5} \cite{H1}
\cite{H6} and the proof of the Verlinde conjecture by the author \cite{H10},
which in turn was based on the constructions and studies of 
genus-zero and genus-one correlation functions in \cite{H6} and 
\cite{H7}. 

The semisimplicity of the module categories 
in the work above simplifies many constructions and 
proofs (although, as the required material listed above indicates,
the semisimple theory is already substantial and highly nontrivial). 
But the theory of semisimple finite tensor categories is far from 
the whole story. 
The study of finite tensor categories in \cite{EO} is in fact motivated
by generalizations of the results in the semisimple case to the
general case. In the case of conformal field theories, 
nonsemisimple generalizations of rational conformal field theories are 
called ``logarithmic conformal field theories.''
Logarithmic operator product expansions were first studied 
by Gurarie \cite{G} and logarithmic
conformal field theory has been developed rapidly in recent years
by both physicists and mathematicians.
It has been conjectured that certain candidates for logarithmic 
conformal field theories should give finite tensor categories. 
But not even a precise formulation of a general conjecture 
has been previously given in the literature. 

In the present paper, we discuss what has been done and 
what still needs to be done for the problem of constructing 
braided finite tensor category structure on a suitable module category
for a vertex operator algebra. 
We shall also identify clearly the main difficulties that we have encountered 
in establishing
these results and the methods that we have developed to overcome them.
We will also present problems that still need to be solved
and discuss three applications. In particular, we give a general 
conjecture on the class of vertex operator algebras for which 
the categories 
of grading-restricted generalized modules have natural structures
of finite tensor categories.

The present paper is organized as follows: We recall briefly 
some basic notions in the theory of tensor categories in 
Section 2. In Section 3, we discuss results 
on the construction of a modular tensor category from 
modules for a vertex operator algebra satisfying certain natural 
positive energy, finiteness and reductivity conditions.
In Section 4, we discuss results on the construction of a braided 
tensor category from grading-restricted generalized modules 
for a vertex operator algebra satisfying certain positive energy
and finiteness conditions (but not necessarily the reductivity
condition). Conjectures and problems are also discussed in this section. 
Applications are 
discussed in Section 5. 

\paragraph{Acknowledgment} This paper is dedicated to Geoffrey 
Mason for his 60th birthday.  The author would like to 
thank J\"{u}rgen Fuchs, Liang Kong, Jim Lepowsky,
Antun Milas and Christoph Schweigert for comments.

\renewcommand{\theequation}{\thesection.\arabic{equation}}
\renewcommand{\thethm}{\thesection.\arabic{thm}}
\setcounter{equation}{0}
\setcounter{thm}{0}

\section{Finite tensor categories}

We first recall the basic notions in the theory of tensor 
categories. 
The purpose is mainly to clarify the terminology because 
different terminologies exist in the literature. 
We shall be sketchy in describing these notions.
See \cite{T2}, \cite{BK} and \cite{EO} for more details. 

A {\it tensor category} is an abelian category with a 
monoidal category structure. A {\it braided tensor category}
is a tensor category with a natural braiding isomorphism
from the tensor product bifunctor to the composition of 
the tensor product bifunctor and the permutation functor on the 
direct product of two copies of the category,
such that the two standard hexagon diagrams are commutative (see \cite{T2}
and \cite{BK}). 
A tensor category with tensor product bifunctor $\boxtimes$ 
and unit object $V$ 
is {\it rigid} if for every object 
$W$ in the category, there are right and left dual objects $W^{*}$
and $^{*}W$ together with morphisms $e_{W}: W^{*}\boxtimes W\to V$,
$i_{W}:V\to W\boxtimes W^{*}$, $e'_{W}: W\boxtimes {}^{*}W\to V$
and $i'_{W}: V\to {}^{*}W\boxtimes W$ 
such that the compositions of the morphisms in the sequence
$$\begin{CD}
W&@>>>
&V\boxtimes W
&@>i_{W}\boxtimes 1_{W}>>
&(W\boxtimes W^{*})\boxtimes 
W&@>>>\\
&@>>>&W\boxtimes (W^{*}\boxtimes 
W)&@>1_{W}\boxtimes e_{W}>>&
W\boxtimes V&@>>>&W
\end{CD}$$
and three similar sequences are equal to the identity
$1_{W}$ on $W$ (see \cite{T2} and \cite{BK}). 
A rigid braided tensor category together with a twist (a natural 
isomorphism from the category to itself) satisfying natural conditions 
(see \cite{T2} and \cite{BK} for the precise conditions)
is called a {\it ribbon category}. 

An object $W$ in an abelian category is {\it simple} if 
any monomorphism to $W$ is either $0$ or an isomorphism.
An abelian category is said to be {\it semisimple} if every 
object is isomorphic to a direct sum of simple objects.
An object $W$ is of {\it finite length} if there exists a 
finite sequence of monomorphisms $0\to W_{n}\to \cdots\to W_{0}=W$
such that the cokernels of these monomorphisms are simple 
objects.

An object $W$ in an abelian category is {\it projective}
if for any objects $W_{1}$ and $W_{2}$, any morphism
$p: W\to W_{2}$ and any epimorphism $q: W_{1}\to W_{2}$, 
there exists a morphism
$\tilde{p}: W\to W_{1}$ such that $q\circ \tilde{p}=p$. 
Let $W$ be an object of the abelian category. A {\it projective cover
of $W$ in the category} is a projective object $U$
and an epimorphism $p: U\to W$ such that 
for any projective object $W_{1}$ and any 
epimorphism $q: W_{1}\to W$, there exists an
epimorphism $\tilde{q}: W_{1}\to U$ such that $p\circ \tilde{q}=q$.

A {\it finite tensor category} is a rigid tensor category 
such that every object is of finite length, every space of 
morphisms is finite-dimensional, there are only finitely many 
inequivalent simple objects, and every simple
object has a projective cover. A {\it braided finite tensor 
category} is a finite tensor category
which is also a braided tensor category.

\begin{expl}
{\rm The category of finite-dimensional modules for a finite group is 
a finite tensor category. But the category of finite-dimensional 
modules for a simple finite-dimensional Lie algebra is not a 
finite tensor category.}
\end{expl}

\section{The semisimple case}

A semisimple ribbon category with 
finitely many inequivalent simple objects $W_{1}, \dots, W_{m}$
and braiding isomorphism $c$ is a {\it modular 
tensor category} if it has the following nondegeneracy 
property:
The $m\times m$ matrix formed by the traces
of the morphisms $c_{W_{i}W_{j}}\circ c_{W_{j}W_{i}}$ 
in the ribbon category for $i, j=1, \dots, m$ is invertible. 
See \cite{T2} and \cite{BK} for studies of modular tensor 
categories. 
In this semisimple case, simple objects are projective 
covers of themselves. 

In 1988, as we mentioned in the introduction, 
Moore and Seiberg \cite{MS1} \cite{MS2}
derived a set of polynomial equations
from an axiom system for a rational conformal field theory. Inspired
by a comment of Witten, they observed an analogy between 
the theory of these polynomial equations and the theory of 
tensor categories. The structures given by these Moore-Seiberg equations
were called modular tensor categories by I. Frenkel.
However,  in the work of Moore and Seiberg, as they commented,
tensor product and other structures were neither formulated nor
constructed mathematically. Later, 
Turaev formulated a notion of modular tensor
category in \cite{T1} and \cite{T2} and gave examples 
of such tensor categories
from representations of quantum groups at roots of unity
based on results obtained by many people on quantum groups and their
representations, especially those in the pioneering work 
\cite{RT1} and \cite{RT2} by Reshetikhin and Turaev 
on the construction of knot and $3$-manifold invariants from 
representations of quantum groups. 
The original structures
given by the Moore-Seiberg equations can be obtained easily from 
modular tensor categories
in this sense and are analogous to 
$6$-$j$ symbols in the representation theory of Lie algebras. This 
new and conceptual formulation of the notion of 
modular tensor category by Turaev
led to the conjecture that a rational conformal field theory 
(if such a structure actually exists)
gives naturally a modular tensor category in this sense of Turaev. 
Moreover, since  
the construction of rational conformal field theories
is  harder than the construction of modular tensor 
categories, a more appropriate 
problem is to construct directly modular tensor categories in the sense of 
Turaev from representations of vertex operator algebras, 
which are substructures of candidates for rational conformal field
theories. In fact, it turns out that a series of results obtained 
by the author and collaborators in the construction of modular tensor categories 
from representations of vertex operator algebras are necessary 
steps (but already sufficient for many applications so far) 
in the author's program of constructing rational conformal field theories
from representations of vertex operator algebras.

In physics, there have been known candidates for rational conformal 
field theories, for example, the Wess-Zumino-Novikov-Witten (WZNW)
models
and the minimal models. Until 2005, 
it was a well-known conjecture 
that the categories of suitable modules for affine 
Lie algebras of positive integral levels and for the Virasoro algebra
of certain central charges and some other categories studied by physicists
are indeed modular tensor categories. 

The first mathematical construction of a rigid braided tensor category 
structure 
from representations of affine Lie algebras was given by Kazhdan and 
Lusztig \cite{KL1}--\cite{KL5}. But these tensor categories are not 
finite and do not correspond to rational conformal field theories. 
Under the assumption that 
the conjectured braided tensor category structure 
on the category of integrable highest weight modules of a positive
integral level for an affine Lie algebra is rigid, Finkelberg \cite{Fi1}
\cite{Fi2} showed that 
this conjectured braided tensor category structure can actually be 
obtained by transporting the corresponding
braided tensor category structure constructed by Kazhdan and Lusztig
to this category. There were also the works of Tsuchiya-Ueno-Yamada
\cite{TUY} and Beilinson-Feigin-Mazur \cite{BFM}, in which the WZNW models 
and minimal models were studied using algebro-geometric methods. 
The author was told by experts
that the results obtained in these works can be used
to construct braided tensor category structures on the 
corresponding module categories for WZNW models and minimal models.
However, the rigidity and the nondegeneracy 
property of these braided tensor categories 
cannot be proved using the results and 
methods in \cite{TUY} and \cite{BFM}. The book \cite{BK} 
gave a construction of the braided
tensor categories for WZNW models but did not give a proof of 
the rigidity. 
So even in the case of 
WZNW and minimal models, the construction of 
the corresponding modular tensor category structures 
was still an unsolved open problem before 2005. 

On the other hand, starting from 1991, Lepowsky and the author
in \cite{tensor0}--\cite{tensor3}, \cite{tensor5} and \cite{H1}
developed a tensor product theory for modules for a vertex operator
algebra satisfying suitable finiteness and reductivity conditions.
In particular, the braided tensor category structures 
for the WZNW and minimal models were constructed in \cite{affine}
and \cite{H2}, respectively, using this general theory. In this 
tensor product theory, the 
hard part is the construction of the associativity isomorphism and 
the proof of the commutativity of the pentagon and hexagon diagrams
(the main coherence properties). 
The method used required both  algebra (especially the method 
of formal variables) and  complex analysis. For example, 
one of the main steps in the construction of the associativity 
isomorphism and the proof of the coherence properties is the 
proof of the convergence of products of intertwining operators. 
The formal variable method is necessary because to prove
the convergence, we have to prove that
the formal series of products of intertwining operators 
satisfy differential equations with formal series as coefficients.
On the other hand, we cannot construct
the associativity isomorphism and prove the coherence properties
without using some delicate  complex analysis.
The complex variable method is necessary and no algebraic method 
can be used to replace it. We not only have to prove the 
convergence, but also have to deal with very subtle issues for 
the convergent series.
For example, if $\sum_{n\in D}a_{n}z^{n}=0$, is it true that 
$a_{n}=0$ for $n\in D$? If $D=\C$, then the answer is no. 
If $D$ is a discrete subset of $\R$, 
then the answer is yes. In the original construction of the 
associativity isomorphism in \cite{H1}, $D$ is 
assumed to be a strictly increasing sequence in $\R$. In particular, 
$D$ in this case is discrete. 

As in the special cases of WZNW and minimal models, 
the rigidity of these braided tensor categories was still a conjecture
before 2005. We now know that the reason why the rigidity 
was so hard is that one needs
the Verlinde conjecture to prove the rigidity. The proof 
of the Verlinde conjecture by the author
in \cite{H10} requires not only the genus-zero
theory (the theory of intertwining operators) 
but also the genus-one theory (the theory of $q$-traces of 
intertwining operators and their modular invariance). 
Note that the statement of rigidity actually involves only the 
genus-zero theory but its proof in \cite{H11} needs the genus-one theory. 
There must be something deep going on here. 

The modular invariance result needed in the proofs of the Verlinde conjecture,
the rigidity and the nondegeneracy property is the (strong) result for 
intertwining operators established in 
\cite{H7}. The modular invariance proved by Zhu \cite{Zhu1}
\cite{Zhu2} is only a 
very special case of this stronger result needed, and 
is far from enough for these purposes. The paper \cite{H7} not only 
established the most general modular invariance result in 
the semisimple case, but also constructed all genus-one correlation 
functions of the corresponding chiral rational conformal 
field theories. After Zhu's modular invariance was proved
in 1990, the 
modular invariance for products or iterates of more than one 
intertwining operator had been an open problem for a long time. 
In the case of products or iterates of at most one intertwining 
operator and any number of vertex operators for modules, 
a straightforward generalization of Zhu's result using the same method
gives the modular invariance (see \cite{M}). 
But for products or iterates of more than one 
intertwining operator, Zhu's method simply does not work. In fact, in 
this general case, even the theory of intertwining operators 
had not been fully developed before 2003. This is 
one of the main reasons that for about 15 years
after 1990, there had been not 
much progress towards the proof of the rigidity and 
the nondegeneracy property. 

This situation changed in 2003 when the author constructed chiral
genus-zero correlation functions using intertwining operators \cite{H6}
and proved the modular invariance of the space of $q$-traces 
of products and iterates 
of intertwining operators \cite{H7}. These results are for a simple vertex 
operator algebra $V$  satisfying the following conditions (for 
basic definitions and terminology in the theory 
of vertex operator algebras, see \cite{FLM},
\cite{FHL}, \cite{LL}, \cite{HLZ2} and \cite{H12}):
\renewcommand{\labelenumi}{\Roman{enumi}}
\begin{enumerate}

\item $V$ is of positive energy ($V_{(0)}=\C\one$ and $V_{(n)}=0$ 
for $n<0$) and the contragredient $V'$, as a $V$-module, is equivalent to $V$.

\item Every $\N$-gradable weak $V$-module is a direct sum of 
irreducible $V$-modules. 
(In fact, the results proved in 
\cite{H12} imply that this condition can be weakened to 
the condition that every grading-restricted 
generalized $V$-module is a direct sum of 
irreducible $V$-modules.)

\item $V$ is $C_{2}$-cofinite.

\end{enumerate}
Using these results the author 
proved the Verlinde conjecture in \cite{H10}
in 2004
and  the rigidity and 
nondegeneracy property for the braided tensor category 
of modules for a vertex operator algebra satisfying Conditions
I--III in \cite{H11} in 2005\footnote{In fact, the author 
proved the rigidity and 
nondegeneracy property before the summer of 2004 and 
discussed the proof in talks in two conferences in 2004. But the
paper \cite{H11} was posted to the arXiv in 2005.}. In particular, we have:

\begin{thm}[\cite{H11}]\label{3.1}
Let $V$ be a simple vertex operator algebra satisfying
Conditions I--III above.
Then the category of $V$-modules has a natural structure of 
modular tensor category.
\end{thm}

\begin{rema}
{\rm From this result and \cite{T2}, we obtain a modular 
functor (including all genus) which in particular gives a 
representation of the mapping class group
of a Riemann surface of any genus. On the other hand, 
for the vertex operator algebra $V$,
chiral correlation functions (or conformal blocks) 
on Riemann surfaces of any genus 
can also be defined directly. These chiral correlation functions
or conformal blocks 
also give a representation of the mapping class group of a 
Riemann surface of any genus. These two representations 
of the mapping class group
of the same Riemann surface are certainly expected to be the same. 
But the proof of this fact needs a construction of the chiral correlation
functions on higher-genus Riemann surfaces
from intertwining operators. This construction is
still an unsolved problem. In fact, the only unsolved part
in this construction is a suitable convergence problem,
which is now also the main unsolved problem 
in the author's program of constructing higher-genus 
chiral rational conformal field theories from 
representations of vertex operator algebras satisfying Conditions
I--III. }
\end{rema}

\section{The general (not necessarily semisimple) case}

The semisimplicity of the categories under consideration, 
as we have discussed 
in the preceding section, simplifies many things.
But it is not natural to study only the semisimple case.
The semisimple case is also not general enough for
further developments and applications.
A satisfying theory should be a theory in the general case (not necessarily 
semisimple) and the theory in the semisimple case would become
a special case of the general theory.

To achieve this, we need to remove Condition II
discussed in the preceding section. In this general case, our 
strategy is still the same as in the semisimple case. 
The first step is to construct braided tensor category structures.
The second step is to prove the rigidity. The last step is 
to formulate and prove a generalization of the nondegeneracy 
property in this general case. 
The first step has been carried out recently by the author in \cite{H12}
using the general logarithmic tensor product theory (in which 
the semisimple theory in \cite{tensor0}--\cite{tensor3}
\cite{tensor5}, \cite{H1} and 
\cite{H6} is indeed a special case) developed
by Lepowsky, Zhang and the author \cite{HLZ1} \cite{HLZ2} and 
a number of results (see \cite{H12} for details) 
in the representation theory of vertex operator
algebras obtained by many people in the past twenty years.
The second and third steps will need modular invariance 
in the general (not necessarily semisimple) case and 
are one of the research projects that the author is finishing. 

In the semisimple case, we know from \cite{H6} and \cite{H11}
that
the $C_{1}$-cofiniteness conditions together with some 
other minor conditions are enough for the construction of
braided tensor category structures 
while the $C_{2}$-cofiniteness conditions are needed
only in the proof of the rigidity and the nondegeneracy 
property. In the general case, for braided tensor category structures,
we also do not need
the stronger $C_{2}$-cofiniteness condition;
$C_{1}$-cofiniteness
conditions together with some other conditions are enough.
More precisely, we consider the following conditions for a vertex operator
algebra $V$:
\renewcommand{\labelenumi}{\arabic{enumi}}
\begin{enumerate}

\item There exists a positive integer $N$ such that the difference between
the lowest weights of any two irreducible $V$-modules is less than $N$ 
and the associative algebra $A_{N}(V)$ (see \cite{DLM} and the explanations
below)
is finite dimensional. 

\item $V$ is $C_{1}$-cofinite in the sense of Li \cite{L}.

\item Irreducible $V$-modules are $\R$-graded and are $C_{1}$-cofinite 
in the sense of \cite{H6} (or quasi-rational in the sense of Nahm
\cite{N}). 

\end{enumerate}

Here are some explanations of the conditions above:
\begin{enumerate}

\item $A_{N}(V)$ is the natural generalization 
of Zhu's algebra (see \cite{Zhu1} and \cite{Zhu2})
by Dong-Li-Mason \cite{DLM}. 

\item  Li's $C_{1}$-cofiniteness condition in \cite{L}
can also be defined for modules,
but it is mainly useful for the algebra $V$. We believe that this should 
be the correct cofiniteness condition on $V$ needed 
for genus-zero theories.

\item The  $C_{1}$-cofiniteness condition for 
modules used in \cite{H6} 
always holds for $V$. It was introduced first by Nahm in \cite{N},
where such modules are called quasi-rational. 
Clearly it cannot be used as a condition
for the algebra. But this 
$C_{1}$-cofiniteness condition for modules
is important for getting differential equations satisfied
by intertwining operators. We believe that this is the correct
cofiniteness condition for modules needed for genus-zero theories. 

\item If $V$ satisfies Conditions I and III in the preceding section, 
it satisfies Conditions 1--3 above. 

\end{enumerate}

The logarithmic tensor product theory 
for suitable categories of generalized modules for vertex operator
algebras developed by Lepowsky, Zhang and the author \cite{HLZ1}
\cite{HLZ2} says that 
if the vertex operator algebra $V$, generalized 
$V$-modules in a suitable category $\mathcal{C}$
and logarithmic intertwining 
operators satisfy certain conditions,
including in particular the assumption that $\mathcal{C}$ is closed
under a candidate for the tensor product bifunctor, then the category 
$\mathcal{C}$ has a natural structure of a braided tensor category.
Using this general theory, the following result is obtained by the author:

\begin{thm}[\cite{H12}]\label{4.1}
Assume that $V$  satisfies
Conditions 1--3 above. Then 
every irreducible $V$-module has a projective cover and
the category of grading-restricted generalized $V$-modules has a 
natural structure of a braided tensor category.
\end{thm}

The logarithmic tensor product theory
in \cite{HLZ1} and \cite{HLZ2} reduces the proof of this 
theorem to the proof that the conditions
needed to use the logarithmic tensor product theory 
are all satisfied. In this case, the proof of
the convergence and extension properties are similar to 
the semisimple case. The hard part is the proof  that 
the tensor product of two grading-restricted generalized $V$-modules
is still a grading-restricted generalized $V$-module, or equivalently,
the existence of the tensor product bifunctor.
The crucial steps in this proof are the proof of 
the existence of projective covers using 
the theory of $A_{N}(V)$-algebras
and the proof of the existence of the tensor product bifunctor
using projective covers.

The theorem above does not say anything about the 
rigidity of the braided tensor category. 
The author believes that in general the rigidity 
is not true for these braided tensor categories. But we have 
the following:

\begin{conj}\label{4.2}
Assume that $V$ is a simple vertex operator algebra 
satisfying Condition I and III in the preceding 
section; in particular, Conditions 1--3 above hold. Then 
the braided tensor category given in Theorem \ref{4.1}
is rigid. 
\end{conj}

This conjecture immediately implies the following:

\begin{cor}\label{4.3}
Assume that $V$ is a simple vertex operator algebra 
satisfying Condition I and III in the preceding 
section. Then the category of grading-restricted generalized 
$V$-modules has a 
natural structure of braided finite tensor category. Moreover,
equipped with a natural twisting, it is a ribbon category. 
\end{cor}

The proof of Conjecture \ref{4.2} is expected to be similar 
to Theorem \ref{3.1} in the semisimple case. In that case,
the rigidity was proved using the Verlinde conjecture which in turn
was proved using the modular invariance and the genus-one 
associativity established in \cite{H7}. So to prove  Conjecture \ref{4.2},
we need first to establish a generalization of the modular invariance
and the genus-one associativity  and then to 
formulate and prove a generalization of the Verlinde conjecture relating 
fusion rules and modular transformations, both in this general case. 
The author has generalized his proofs of the modular invariance and 
genus-one associativity in the semisimple case to this general case 
and a generalization of the Verlinde conjecture
is expected to be a consequence of these results. 
As in the semisimple case, the rigidity will be a consequence of 
these results. 
We know that in the semisimple case the 
nondegeneracy property also follows easily from the Verlinde 
conjecture. Therefore we expect that a proof 
of a suitable nondegeneracy property in this case 
can  be obtained 
based on these results.

\section{Applications}

The tensor category structures discussed in the preceding 
two sections have many applications. In fact, these tensor 
category structures had been conjectured to exist for many years.
Many results were obtained by physicists and mathematicians
based on this postulated existence. These  results 
are of their own importance and 
interest in different areas of mathematics and physics. 
Here we choose to discuss three of these applications 
in this section.

\subsection{Open-closed conformal field theories}

The first application is on algebras in modular tensor categories
and open-closed conformal field theories.
Various notions of algebra, including associative algebra,
commutative associative algebra, Frobenius algebra and so on, can be
defined in a braided tensor category. In the case that the modular tensor
category is constructed from the category of modules for a vertex
operator algebra in Theorem \ref{3.1}, these algebras 
are equivalent to substructures
of open-closed conformal field theories. 
Therefore we can apply our theory of these modular tensor categories
to the study of 
open-closed conformal field theories. 

Comparing tensor categories constructed from representations 
of vertex operator algebras with tensor categories
of vector spaces,  we can see that vertex operator algebras in Theorems
\ref{3.1} and \ref{4.1} in fact
play a role analogous to the coefficient fields of vector spaces. 
Given a field, we have a symmetric tensor category of vector spaces
over this field. The theory of all types of algebras is
based on this symmetric tensor category structure. 
In our case, given a vertex operator algebra satisfying the conditions in 
Theorems \ref{3.1} or \ref{4.1}, we have a braided tensor category.
The general theory of algebras in braided tensor categories can now be applied 
to study algebras in this particular braided tensor category. 

Under suitable assumptions, including, in particular,
the existence of a modular tensor category structure on 
the category of modules for a vertex operator algebra,
in a series of papers, Felder, Fr\"{o}hlich,
Fuchs, Schweigert, Fjelstad and Runkel \cite{FFFS} \cite{FFRS1}
\cite{FFRS2} \cite{FRS1} \cite{FRS2}
developed an approach to open-closed
conformal field theories using algebras in modular 
tensor categories and 
three-dimensional topological field theories constructed from 
such categories. 
Because of Theorem \ref{3.1},
the results on modular tensor categories and on algebras in 
these tensor categories in these papers are indeed equivalent to 
results in open-closed conformal field theory.  
In this approach, one starts with a modular tensor
category and a symmetric Frobenius algebra in this category
and constructs correlators for the corresponding
open-closed conformal field theory from the category and
the algebra. 

There is another approach to open-closed conformal field theories 
developed by the author \cite{H2}--\cite{H5}, by Kong and the author 
\cite{HK1}--\cite{HK3} and by Kong \cite{Ko1}--\cite{Ko3}
using directly the representation theory of vertex operator algebras.
In this approach, 
one starts from a vertex operator algebra satisfying suitable
conditions and constructs correlation functions of the corresponding
open-closed conformal field theory from representations of the vertex
operator algebra. 
Using the tensor product theory developed by Lepowsky and the 
author \cite{tensor0}--\cite{tensor3}
\cite{tensor5} \cite{H1}
\cite{H6} and the author's construction of the modular tensor category 
structures in \cite{H11}, 
Kirillov, Lepowsky and the author \cite{HKL},  Kong and the 
author \cite{HK1} and Kong \cite{Ko1} 
established the equivalence of suitable algebras 
in suitable modular tensor categories with suitable vertex operator
algebras, open-string vertex operator algebras or full field algebras, 
respectively.
Using all these results, Kong \cite{Ko3} introduced 
what he called Cardy algebras which are conjectured to 
be equivalent to open-closed conformal field theories.

In a recent paper \cite{KR}, Kong and Runkel studied the relations between 
these two approaches discussed above and unified them in a single
framework for
open-closed conformal field theories.

The results discussed above on algebras in braided tensor 
categories are all given in the semisimple case, since the 
corresponding open-closed  conformal field theories 
are rational. But many results can be easily generalized to the 
general case. It will be interesting to see how much of 
the theory in the semisimple case can be generalized to the 
not-necessarily-semisimple case. 

\subsection{Triplet $W$-algebras}

The second application is on a conjectured equivalence
between the braided finite tensor category 
of grading-restricted generalized modules for a triplet $W$-algebra
and the braided finite tensor category of suitable 
modules for a restricted quantum group. 

The triplet $\mathcal{W}$-algebras
of central charge $1-6\frac{(p-1)^{2}}{p}$ were
introduced first by Kausch \cite{K1} and have been studied extensively
by Flohr \cite{F1} \cite{F2}, 
Gaberdiel-Kausch \cite{GK1} \cite{GK2}, Kausch \cite{K2}, 
Fuchs-Hwang-Semikhatov-Tipunin \cite{FHST}, Abe \cite{A},
Feigin-Ga{\u\i}nutdinov-Semikhatov-Tipunin \cite{FGST1} \cite{FGST2}
\cite{FGST3}, Carqueville-Flohr \cite{CF}, Flohr-Gaberdiel 
\cite{FG}, Fuchs \cite{Fu}, Adamovi\'{c}-Milas \cite{AM1} \cite{AM2},
Flohr-Grabow-Koehn \cite{FGK}, 
Flohr-Knuth \cite{FK} and Gaberdiel-Runkel \cite{GR1} \cite{GR2}.
Based on the results of Feigin-Ga{\u\i}nutdinov-Semikhatov-Tipunin obtained 
\cite{FGST1}  
and of Fuchs-Hwang-Semikhatov-Tipunin   \cite{FHST}, 
Feigin, Ga{\u\i}nutdinov, Semikhatov and Tipunin 
conjectured \cite{FGST2} the equivalence mentioned above
and proved the conjecture 
in the simplest $p=2$ case. But their formulation of the 
conjecture also includes
the statement that the categories of modules for the 
triplet $\mathcal{W}$-algebras considered in their paper
are indeed braided tensor categories.

The triplet $\mathcal{W}$-algebras are vertex operator algebras
satisfying the positive energy condition
(Conditions I) and the $C_{2}$-cofiniteness condition (Condition III)
but not Condition II.
The $C_{2}$-cofiniteness condition was proved by Abe \cite{A}
in the simplest $p=2$ case
and by Carqueville-Flohr \cite{CF} and Adamovi\'{c}-Milas
\cite{AM2} in the general case. 
Condition II was proved to be not satisfied by 
these vertex operator algebras by Abe  \cite{A} in the 
simplest  $p=2$ case
and by Fuchs-Hwang-Semikhatov-Tipunin \cite{FHST} and
Adamovi\'{c}-Milas \cite{AM2} in the general case. 
By Theorem \ref{4.1}, the category 
of grading-restricted generalized modules for 
a triplet $W$-algebra has a natural braided tensor category
structure, not just a quasi-tensor category structure. 
Assuming that Conjecture \ref{4.2} is true for 
a triplet $W$-algebra, we see that by 
Theorem \ref{4.1} and Corollary \ref{4.3},
the category of grading-restricted generalized modules 
for a triplet $W$-algebra is a braided finite tensor category. 
Then the conjecture of Feigin, Gainutdinov, Semikhatov and Tipunin  
is now purely about the equivalence between 
the braided finite tensor category 
of grading-restricted generalized modules for a triplet $W$-algebra
and
the braided finite tensor category of suitable modules for the corresponding
restricted quantum group. 
We expect that the logarithmic 
tensor product theory 
developed in \cite{HLZ1}, \cite{HLZ2} and \cite{H12} 
will be useful in proving and understanding 
this conjecture.

\subsection{Knot and $3$-manifold invariants}

Finally we discuss the application to the construction and 
study of $3$-dimensional topological field theories 
and knot and $3$-manifold invariants.
In \cite{T2}, Turaev constructed $3$-dimensional topological field theories 
and knot and $3$-manifold invariants from modular tensor categories. 
Combining this result of Turaev with Theorem \ref{3.1}, 
we immediately obtain a $3$-dimensional topological field theory,
and a knot and $3$-manifold invariant from a simple vertex operator
algebra satisfying the condition I--III needed in Theorem \ref{3.1}. 

In \cite{He}, Hennings constructed a topological invariant of $3$-manifolds
 from  quantum groups in a manner 
similar to the Witten-Reshetikhin-Turaev invariant \cite{W} \cite{RT2}.
But in this construction, 
instead of working with a semisimple part of the category of 
the representations of a quantum group, Hennings worked directly with 
the nonsemisimple theory. See also the refinement by Kauffman and Radford
in \cite{KaR}.
In \cite{KL}, Kerler and Lyubashenko constructed a $3$-dimensional
extended topological field theory from a modular bounded 
abelian ribbon category, which is a nonsemisimple generalization of 
a modular tensor category. In \cite{Ke}, it was shown that 
underlying the Hennings invariant is exactly the nonsemisimple $3$-dimensional
extended topological field theories constructed in \cite{KL} (see also 
\cite{Ke1}). 

If Conjecture \ref{4.2} and thus 
Corollary \ref{4.3} is proved and the nondegeneracy property 
is formulated and proved, we might obtain a modular bounded 
abelian ribbon category in the sense of Kerler and Lyubashenko
or a similar structure. It is reasonable to conjecture that 
we should be able to obtain a $3$-dimensional 
extended  topological field theory 
or some other natural generalization of a semisimple $3$-dimensional
topological field theory. It will be interesting to see whether 
we will be able to obtain new knot and $3$-manifold invariants
in this way.

\noindent {\small \sc Department of Mathematics, Rutgers University,
110 Frelinghuysen Rd., Piscataway, NJ 08854-8019}

\noindent {\em E-mail address}: yzhuang@math.rutgers.edu


\begin{thebibliography}{KWAK2}

\bibitem[A]{A} 
T. Abe,
A $\mathbb{Z}\sb 2$-orbifold model of the symplectic fermionic vertex
operator superalgebra, {\it Math. Z.} {\bf 255} (2007), 755--792.

\bibitem[AM1]{AM1} D. Adamovi\'c and A. Milas, Logarithmic 
intertwining operators
and $\mathcal{W}(2,2p-1)$-algebras, {\it J. Math. Phys.}
{\bf 48}, 073503 (2007).

\bibitem[AM2]{AM2} D. Adamovi\'c and A. Milas, 
On the triplet vertex algebra $\mathcal{W}(p)$, 
{\it Adv. in Math.} {\bf 217}  (2008), 2664-2699.



\bibitem[BK]{BK}
B. Bakalov and A. Kirillov, Jr., {\em Lectures 
on tensor categories and modular functors}, 
University Lecture Series, Vol. 21,
Amer. Math. Soc., Providence, RI, 2001.

\bibitem[BFM]{BFM} 
A. Beilinson, B. Feigin and B. Mazur, Introduction to algebraic field
theory on curves, preprint, 1991 (provided by A. Beilinson, 1996).

\bibitem[CF]{CF} 
N. Carqueville and M. Flohr, Nonmeromorphic operator
product expansion and $C_2$-cofiniteness for a family of
$\cal{W}$-algebras, {\it J.Phys.} {\bf A39} (2006), 951--966.


\bibitem[DLM]{DLM} 
C. Dong, H. Li and G. Mason, Vertex operator algebras and 
associative algebras, {\it J. Algebra} {\bf 206} (1998), 67--96.



\bibitem[EO]{EO} 
P. Etingof and V. Ostrik, Finite tensor categories,
{\it Moscow Math. J.} {\bf 4} (2004), 627--654.

\bibitem [FGST1]{FGST1} 
B. L. Feigin, A. M. Ga\u\i nutdinov, 
A. M. Semikhatov, and I. Yu Tipunin, I,
The Kazhdan-Lusztig correspondence for the representation category
of the triplet $W$-algebra in logarithmic conformal field theories
(Russian), {\it Teoret. Mat. Fiz.} {\bf 148} (2006), no. 3, 398--427.

\bibitem [FGST2]{FGST2} 
B. L. Feigin, A. M. Ga\u\i nutdinov, 
A. M. Semikhatov, and I. Yu Tipunin,
Logarithmic extensions of minimal models: characters and modular
transformations, {\it Nucl. Phys.} B {\bf 757} (2006), 303--343.

\bibitem [FGST3]{FGST3} 
B. L. Feigin, A. M. Ga\u\i nutdinov, 
A. M. Semikhatov, and I. Yu Tipunin,
Modular group representations and fusion in logarithmic conformal
field theories and in the quantum group center, {\it Comm. Math.
Phys.} {\bf 265} (2006), 47--93.

\bibitem[FFFS]{FFFS}
G. Felder, J. Fr\"{o}hlich, J. Fuchs and C. Schweigert, 
Correlation functions and boundary
conditions in rational conformal field theory and three-dimensional topology, 
{\it Compositio
Math.} {\bf 131} (2002) 189--237.

\bibitem[Fi1]{Fi1}
M. Finkelberg, Fusion categories, Ph.D. thesis, Harvard University,
1993.

\bibitem[Fi2]{Fi2}
M. Finkelberg, An equivalence of fusion categories, 
{\em Geom. Funct. Anal.} {\bf 6} (1996), 
249--267.

\bibitem[FFRS1]{FFRS1}
J. Fjelstad, J. Fuchs, I. Runkel, and C. Schweigert,
Uniqueness of open/closed rational CFT with given algebra of open states,
{\it Adv. Theor. Math. Phys.} {\bf 12} (2008), 1283-1375. 

\bibitem[FFRS2]{FFRS2}
J. Fjelstad, J. Fuchs, I. Runkel, and C. Schweigert,
TFT construction of RCFT correlators V: Proof of modular invariance 
and factorisation, {\it Theo. Appl. Categories} {\bf 16} (2006) 342-433.

\bibitem[Fl1]{F1}
M. Flohr, On modular invariant partition functions of 
conformal field theories with logarithmic operators, 
{\it Int. J. Mod. Phys.} {\bf A11} (1996), 4147--4172.

\bibitem[Fl2]{F2}
M. Flohr, On fusion rules in logarithmic conformal field theories,
{\it Int. J. Mod. Phys.} {\bf A12} (1996), 1943--1958.

\bibitem[FG]{FG}
M. Flohr and M. R. Gaberdiel, Logarithmic torus 
amplitudes, {\it J. Phys.} {\bf A39} (2006), 1955--1968.

\bibitem[FK]{FK}
M. Flohr and H. Knuth, On Verlinde-Like formulas in $c_{p,1}$
logarithmic conformal field theories, to appear; arXiv:0705.0545.

\bibitem[FGK]{FGK}
M. Flohr, C. Grabow and M. Koehn,
Fermionic Expressions for the characters of $c(p,1)$
logarithmic conformal field theories, {\it Nucl. Phys.} {\bf B768} (2007), 
263--276.


\bibitem[FHL]{FHL} I.B. Frenkel, Y.-Z. Huang and J. Lepowsky, On
axiomatic approaches to vertex operator algebras and modules,
preprint, 1989; Memoirs Amer. Math. Soc., Vol. 104, Number 494,
American Math. Soc.  Providence, 1993.



\bibitem[FLM]{FLM} 
I.~B. Frenkel, J.~Lepowsky and A.~Meurman, {\it Vertex
Operator Algebras and the Monster}, Pure and Appl. Math., Vol. 134,
Academic Press, Boston, 1988.

\bibitem[Fu]{Fu} 
J. Fuchs, On nonsemisimple fusion rules and tensor categories, in:
{\it Lie algebras, vertex operator algebras and their applications, 
Proceedings of a conference in honor of James Lepowsky and Robert 
Wilson, 2005}, ed. Y.-Z. Huang and  K. Misra, 
Contemporary Mathematics, Vol. 442, Amer. Math. Soc., Providence, 2007.



\bibitem[FHST]{FHST} J. Fuchs, S. Hwang, A.M. Semikhatov and I. Yu. Tipunin,
Nonsemisimple Fusion Algebras and the Verlinde Formula, {\it Comm.
Math. Phys.} {\bf 247} (2004), no. 3, 713--742.



\bibitem[FRS1]{FRS1}
J. Fuchs, I. Runkel and C. Schweigert, 
TFT construction of RCFT correlators. I: Partition
functions, {\it Nucl. Phys.} {\bf B646} (2002) 353--497.

\bibitem[FRS2]{FRS2}
J. Fuchs, I. Runkel and C. Schweigert, 
TFT construction of RCFT correlators IV: 
Structure constants and correlation functions,
{\it Nucl.Phys.}  {\bf B715} (2005) 539-638.




\bibitem[G]{G} 
V. Gurarie, Logarithmic operators in conformal field
theory, {\it Nucl. Phys.} {\bf B410} (1993), 535--549.

\bibitem[GK1]{GK1}
M. R. Gaberdiel and H. G. Kausch, Indecomposable fusion products,
{\it Nucl. Phys.} {\bf B477} (1996), 298--318.

\bibitem[GK2]{GK2}
M. R. Gaberdiel and H. G. Kausch, A rational logarithmic conformal 
field theory, 
{\it Phys. Lett.} {\bf B386} (1996), 131--137.

\bibitem[GR1]{GR1}
M. R. Gaberdiel and I. Runkel, The logarithmic triplet theory with boundary,
{\it J.Phys.} {\bf A39} (2006), 14745-14780.

\bibitem[GR1]{GR2}
M. R. Gaberdiel and I. Runkel, From boundary to bulk in logarithmic CFT,
{\it J. Phys.} {\bf A41} (2008), 075402.

\bibitem[He]{He}
M. Hennings, Invariants of links and 3-manifolds obtained from 
Hopf algebras, {\it J. London Math. Soc.
(2)} {\bf 54} (1996),  594--624.

\bibitem[H1]{H1}
Y.-Z. Huang, A theory of tensor products for module categories for a
vertex operator algebra, IV, {\it J. Pure Appl. Alg.} 100 (1995)
173--216.

\bibitem[H2]{H1.5}
Y.-Z. Huang, Virasoro vertex operator algebras,
(nonmeromorphic) operator product expansion and the tensor product
theory, {\em J. Alg.} {\bf 182} (1996), 201--234.

\bibitem[H3]{H2} Y.-Z. Huang, Intertwining operator algebras,
genus-zero modular functors and genus-zero conformal field theories,
in: {\it Operads: Proceedings of Renaissance Conferences}, ed. J.-L. Loday,
J. Stasheff, and A. A. Voronov, Contemporary Math., Vol. 202,
Amer. Math. Soc., Providence, 1997, 335--355.

\bibitem[H4]{H3} Y.-Z. Huang, {\it Two-dimensional conformal geometry and
vertex operator algebras}, Progress in Mathematics, Vol. 148,
Birkh\"{a}user, Boston, 1997.

\bibitem[H5]{H4} Y.-Z. Huang, Genus-zero modular functors and
intertwining operator algebras, {\it Internat. J. Math.} 9 (1998), 
845--863.

\bibitem[H6]{H5} Y.-Z. Huang, Riemann surfaces with boundaries and
the theory of vertex operator algebras, in: {\it Vertex Operator
Algebras in Mathematics and Physics}, ed. S. Berman, Y. Billig,
Y.-Z. Huang and J. Lepowsky, Fields Institute Communications, Vol. 39,
Amer. Math. Soc., Providence, 2003, 109--125.


\bibitem[H7]{H6} Y.-Z. Huang, Differential equations and
intertwining operators, {\it Comm. Contemp. Math.} {\bf 7} (2005),
375--400.

\bibitem[H8]{H7} Y.-Z. Huang, Differential equations, duality and
modular invariance, {\it Comm. Contemp. Math.} {\bf 7} (2005), 649--706.

\bibitem[H9]{H8} 
Y.-Z. Huang,  Vertex operator algebras, the Verlinde conjecture 
and modular tensor categories, {\it Proc. Natl. Acad. Sci. USA}
{\bf 102} (2005), 5352--5356. 

\bibitem[H10]{H9} 
Y.-Z. Huang,  Vertex operator algebras, fusion rules and 
modular transformations, in: {\it Non-commutative Geometry and 
Representation Theory in Mathematical Physics}, ed. J. Fuchs, 
J. Mickelsson, G. Rozenblioum and A. Stolin, 
Contemporary Math. Vol. 391, Amer. Math. Soc., Providence, 2005, 
135--148. 

\bibitem[H11]{H10} 
Y.-Z. Huang, Vertex operator algebras and the Verlinde conjecture, 
{\it Comm. Contemp. Math.} {\bf 10} (2008), 103-154. 

\bibitem[H12]{H11} 
Y.-Z. Huang,  Rigidity and modularity of vertex tensor categories, 
{\it Comm. Contemp. Math.} {\bf 10} (2008), 871--911. 

\bibitem[H13]{H12} 
Y.-Z. Huang, Cofiniteness conditions, projective covers and the logarithmic 
tensor product theory, {\it J. Pure Appl. Alg.} {\bf 213} (2009), 458--475. 

\bibitem[HKL]{HKL}
Y.-Z. Huang, A. Kirillov, Jr.
and J. Lepowsky, Braided tensor categories and extensions of vertex operator
algebras, to appear.

\bibitem[HK1]{HK1} 
Y.-Z. Huang and L. Kong,
Open-string vertex algebras, tensor categories and operads
(Yi-Zhi Huang and L. Kong), {\it Comm. Math. Phys.} {\bf 250} (2004), 
433--471. 

\bibitem[HK2]{HK2} 
Y.-Z. Huang and L. Kong, Full field algebras,
{\it Comm. Math. Phys.}{\bf  272}  (2007), 345--396.

\bibitem[HK3]{HK3} 
Y.-Z. Huang and L. Kong, Modular invariance for conformal 
full field algebras, {\it Trans. Amer. Math. Soc.}, to appear;  
arXiv:math/0609570.



\bibitem[HL1]{tensor0}
Y.-Z. Huang and J. Lepowsky, Toward a
theory of tensor products for representations of a vertex operator
algebra, in: {\it Proc. 20th International Conference on Differential
Geometric Methods in Theoretical Physics, New York, 1991},
ed. S. Catto and A. Rocha, World Scientific, Singapore, 1992, Vol. 1,
344--354.

\bibitem[HL2]{tensor-desc}
Y.-Z. Huang and J. Lepowsky, Tensor products of modules for a vertex
operator algebras and vertex tensor categories, in:
     {\it Lie Theory and Geometry,
in honor of Bertram Kostant,}
ed. R. Brylinski, J.-L. Brylinski, V. Guillemin, V. Kac,
Birkh\"{a}user, Boston, 1994, 349--383.

\bibitem[HL3]{tensor1}
Y.-Z. Huang and J. Lepowsky, A theory of tensor products for module
categories for a vertex operator algebra, I, {\it Selecta Mathematica
(New Series)} {\bf 1} (1995), 699--756.

\bibitem[HL4]{tensor2}
Y.-Z. Huang and J. Lepowsky, A theory of tensor products for module
categories for a vertex operator algebra, II, {\it Selecta Mathematica
(New Series)} {\bf 1} (1995), 757--786.

\bibitem[HL5]{tensor3}
Y.-Z. Huang and J. Lepowsky, A theory of tensor
products for module categories for a vertex operator algebra, III,
{\it J. Pure Appl. Alg.} {\bf 100} (1995) 141--171.

\bibitem[HL6]{affine}
Y.-Z. Huang and J. Lepowsky, Intertwining operator algebras and vertex
tensor categories for
affine Lie algebras,  {\it Duke Math. J.} {\bf 99} (1999), 113--134.

\bibitem[HL7]{tensor5}
Y.-Z. Huang and J. Lepowsky, A theory of tensor products for module
categories for a vertex operator algebra, V, to appear.

\bibitem[HLZ1]{HLZ1}
Y.-Z. Huang, J. Lepowsky and L.Zhang, A logarithmic generalization of
tensor product theory for modules for a vertex operator algebra,
{\it Internat. J. Math.} {\bf 17} (2006), 975--1012.

\bibitem[HLLZ]{HLLZ}
Y.-Z. Huang, J. Lepowsky, H. Li and L. Zhang, On the concepts of
intertwining operator and tensor product module in vertex operator
algebra theory, {\it J. Pure Appl. Algebra} {\bf 204} (2006),
507--535.

\bibitem[HLZ2]{HLZ2}
Y.-Z. Huang, J. Lepowsky and L.Zhang, 
Logarithmic tensor product theory for generalized modules for
a conformal vertex algebra, to appear; arXiv:0710.2687.

\bibitem[Ka1]{K1} H. G. Kausch, Extended conformal algebras
generated by multiplet of primary fields, {\it Phys. Lett.} {\bf 259}
B (1991), 448--455.

\bibitem[Ka2]{K2} H. G. Kausch, Symplectic fermions, {\it Nucl. Phys.} B {\bf
583} (2000), 513--541.


\bibitem[KaR]{KaR}
L. H. Kauffman and D. E. Radford,
Invariants of 3-manifolds derived from finite-dimensional Hopf algebras,
{\it J. Knot Theory Ramifications} {\bf 4} (1995), 131--162.

\bibitem[KaL1]{KL1}
D. Kazhdan and G. Lusztig,
Affine Lie algebras and quantum groups,
{\it Duke Math. J., IMRN} {\bf 2} (1991), 21--29.

\bibitem[KaL2]{KL2}
D. Kazhdan and G. Lusztig,
Tensor structures arising {from} affine Lie algebras, I,
{\it J. Amer. Math. Soc.} {\bf 6} (1993), 905--947.

\bibitem[KaL3]{KL3}
D. Kazhdan and G. Lusztig,
Tensor structures arising {from} affine Lie algebras, II,
{\it J. Amer. Math. Soc.} {\bf 6} (1993), 949--1011.

\bibitem[KaL4]{KL4}
D. Kazhdan and G. Lusztig,
Tensor structures arising {from} affine Lie algebras, III, {\it J.
Amer. Math. Soc.} {\bf 7} (1994), 335--381.

\bibitem[KaL5]{KL5}
D. Kazhdan and G. Lusztig,
Tensor structures arising {from} affine Lie algebras, IV,
{\it J. Amer. Math. Soc.} {\bf 7} (1994), 383--453.

\bibitem[Ke1]{Ke1}
T. Kerler, Genealogy of nonperturbative quantum-invariants of 
$3$-manifolds: The surgical family,
in "Geometry and Physics," {\it Lecture Notes in Pure and Applied Physics},
Vol.  184,
Marcel Dekker, 1997,
503-547.

\bibitem[Ke2]{Ke}
T. Kerler, 
Homology TQFT's and the Alexander-Reidemeister 
invariant of $3$-manifolds via Hopf algebras and skein theory,
{\it Canad. J. Math.}  {\bf 55} (2003), 766--821. 

\bibitem[KeL]{KL}
T. Kerler and V.  Lyubashenko, 
{\it Non-semisimple topological quantum field theories for 
$3$-manifolds with corners},
Lecture Notes in Mathematics, Vol. 1765. Springer-Verlag, Berlin, 2001.




\bibitem[Ko1]{Ko1}
L. Kong, Full field algebras, operads and tensor categories,
{\it Adv. in Math.} {\bf 213} (2007), 271--340.

\bibitem[Ko2]{Ko2}
L. Kong,  Open-closed field algebras, 
{\it Comm. Math. Phys.} {\bf 280} (2008), 207-261.

\bibitem[Ko3]{Ko3}
L. Kong, Cardy condition for open-closed field algebras, 	
{\it Comm. Math. Phys.} {\bf 283} (2008), 25-92.

\bibitem[KoR]{KR}
L. Kong and I. Runkel, Cardy algebras and sewing constraints, I, 
to appear; arXiv:0807.3356.


\bibitem[Le]{Le} 
J. Lepowsky, From the representation theory of
vertex operator algebras to modular tensor categories in conformal
field theory, commentary on Y.-Z. Huang's PNAS article ``Vertex
operator algebras, the Verlinde conjecture and modular tensor
categories'', {\it Proc. Nat. Acad. Sci. USA} {\bf 102} (2005),
5304--5305.

\bibitem[LL]{LL} 
J. Lepowsky and H. Li, {\it Introduction to Vertex
Operator Algebras and Their Representations}, Progress in Math.,
Birkh\"auser, Boston, 2003.

\bibitem[Li]{L}
H. Li, Some finiteness properties of regular vertex operator algebras,
{\it J. Alg.} {\bf 212} (1999), 495--514.

\bibitem[M]{M}
M. Miyamoto, Intertwining operators and modular invariance, to appear,
math.QA/0010180.


\bibitem[MS1]{MS1}
G.~Moore and N.~Seiberg,
Polynomial equations for rational conformal field theories,
{\em Phys. Lett.} {\bf B 212} (1988), 451--460.

\bibitem[MS2]{MS2}
G.~Moore and N.~Seiberg,
Classical and quantum conformal field theory,
{\em Comm. Math. Phys.} {\bf 123} (1989), 177--254.

\bibitem[N]{N}
W. Nahm, Quasi-rational fusion products,
{\it Int. J. Mod. Phys.} {\bf B8} (1994), 3693--3702.

\bibitem[RT1]{RT1}
N. Reshetikhin and V. Turaev, Ribbon graphs and their invariants derived 
from quantum groups, {\em Comm. Math. Phys.} {\bf 127} (1990), 1-26.

\bibitem[RT2]{RT2}
N. Reshetikhin and V. Turaev, Invariants of $3$-manifolds 
via link polynomials and quantum groups, {\it Invent. Math.}
{\bf 103} (1991), 547--598.

\bibitem[TUY]{TUY}
A. Tsuchiya, K. Ueno and Y. Yamada, Conformal field theory on universal 
family of stable curves with gauge symmetries, in: {\em Advanced Studies in
Pure Math.}, Vol. 19, Kinokuniya Company Ltd.,
Tokyo, 1989, 459--566.

\bibitem[T1]{T1}
V.  Turaev, Modular categories and $3$-manifold invariants,
{\it Int. J. Mod. Phys.} {\bf B6} (1992), 1807--1824.

\bibitem[T2]{T2}
V. Turaev, {\em Quantum invariants of knots and $3$-manifolds},
de Gruyter Studies in Math., Vol. 18, 
Walter de Gruyter, Berlin, 1994.

\bibitem[W]{W}
E. Witten, Quantum field theory and the Jones polynomial,
{\em Comm. Math. Phys.} {\bf 121} (1989), 351--399.

\bibitem[Zhu1]{Zhu1}
Y. Zhu, Vertex operators, elliptic functions and
modular forms, Ph.D. thesis, Yale University, 1990.

\bibitem[Zhu2]{Zhu2}
Y. Zhu, Modular invariance of characters of vertex operator algebras,
{\it J.
Amer. Math. Soc.} {\bf 9} (1996), 237--307.
\end{thebibliography}
\end{document}